\documentclass[preprint,12pt,3p]{elsarticle}

\usepackage{amssymb}
\usepackage{amsmath}
\usepackage{amsthm}
\newtheorem{definition}{Definition}

\usepackage[utf8]{inputenc}

\usepackage{amsthm}
\theoremstyle{plain} 
\newtheorem{theorem}{Theorem}
\newtheorem{lemma}[theorem]{Lemma}
\newtheorem{corollary}[theorem]{Corollary}

\usepackage{relsize}
\usepackage{bm}
\usepackage{subcaption}
\usepackage{todonotes}
\usepackage{booktabs}
\usepackage{float}

\usepackage{titlesec}
\titleclass{\subsubsubsection}{straight}[\subsection]

\newcounter{subsubsubsection}[subsubsection]
\renewcommand\thesubsubsubsection{\thesubsubsection.\arabic{subsubsubsection}}

\titleformat{\subsubsubsection}
  {\normalfont\normalsize\itshape}{\thesubsubsubsection}{1em}{}
\titlespacing*{\subsubsubsection}{0pt}{3.25ex plus 1ex minus .2ex}{2pt}

\newcommand{\bigsum}[3]{\mathlarger{\mathlarger{\sum}}_{#2}^{#1}{#3}}
\usepackage[colorlinks,plainpages=false]{hyperref}
\hypersetup{
  pdftitle={Statistical Estimation of higher Dedekind Numbers},
  pdfauthor={Alex Fihman, Lennart Van Hirtum, Christian Plessl}
}
\makeatletter
\providecommand*{\toclevel@subsubsubsection}{4}
\makeatother

\journal{Computational Algebra}

\begin{document}

\begin{frontmatter}



\title{Statistical Estimation of higher Dedekind Numbers} 


\author[1]{Alex Fihman\fnref{equal}}
\author[2]{Lennart Van Hirtum\fnref{equal}}
\author[3]{Christian Plessl}
\fntext[equal]{Both authors contributed equally to this work.}
\affiliation[1]{
  organization={Independent Scholar},
  email={alexfihman@gmail.com}
}
\affiliation[2]{
  organization={Paderborn University, Paderborn Center for Parallel Computing},
  city={Paderborn},
  postcode={33098},
  state={North-Rhein Westfalen, Germany},
  email={lennartv@mail.uni-paderborn.de}
}
\affiliation[3]{
  organization={Paderborn University, Paderborn Center for Parallel Computing},
  city={Paderborn},
  postcode={33098},
  state={North-Rhein Westfalen, Germany},
  email={christian.plessl@uni-paderborn.de}
}


\begin{abstract}
We provide highly accurate estimations of the 10th through 15th Dedekind Numbers, to a precision of 4 digits for $D(10)$, to 2 digits for $D(15)$. These estimates were obtained using three methods, including pair matching on large quantities of 9-dimensional monotone Boolean functions for $D(10)$, Reference Subsets for $D(10)$, $D(11)$, and $D(12)$. And our best method "Weight Layer Branching" which provided accurate estimates for all $D(10)$ through $D(15)$, strongly improving on the previous best known estimates by Korshunov and Tian-Shun Chen et al. 
\end{abstract}



\begin{keyword}

Dedekind numbers \sep monotone Boolean functions \sep antichains \sep Monte Carlo methods \sep MCMC \sep combinatorial enumeration \sep computational algebra



\MSC[2020] 06A07 \sep 05A15 \sep 68R05 \sep 68W20 \sep 65C05
\end{keyword}

\end{frontmatter}



\section{Introduction}
Dedekind numbers $D(n)$ count the Monotone Boolean Functions (MBFs) in $n$ variables, or equivalently the antichains of the Boolean lattice $\{0,1\}^n$. They grow doubly exponentially: $\log D(n) \sim \binom{n}{\lfloor n/2 \rfloor}$, so each successive value is astronomically larger than the last. Table~\ref{tab:dedekindNumbers} shows all currently known exact values; $M(n)$ denotes the set of all MBFs in $n$ dimensions, so $D(n) = |M(n)|$. The 32-year gap between $D(8)$ (Wiedemann, 1991 \cite{wiedemannDedekind8}) and $D(9)$ (Van Hirtum et al. (2023) \cite{vanhirtum2023computation} \& Jäkel (2023) \cite{jakel2023computation}) reflects not a lack of effort but the fundamental difficulty of the problem: exact computation requires enumerating a combinatorial structure whose size grows faster than any fixed exponential in $n$.

\begin{table}[H]
    \centering
    \begin{tabular}{c|l|c}
        $D(0)$ & 2 & Dedekind (1897) \\
        $D(1)$ & 3 & Dedekind (1897) \\
        $D(2)$ & 6 & Dedekind (1897) \\
        $D(3)$ & 20 & Dedekind (1897)\\
        $D(4)$ & 168 & Dedekind (1897) \\
        $D(5)$ & 7581 & Church (1940) \\
        $D(6)$ & 7828354 & Ward (1946) \\
        $D(7)$ & 2414682040998 & Church (1965) \\
        $D(8)$ & 56130437228687557907788 & Wiedemann (1991) \cite{wiedemannDedekind8} \\
        $D(9)$ & 286386577668298411128469151667598498812366 & \begin{tabular}[t]{@{}c@{}}
L. Van Hirtum (2023) \cite{vanhirtum2023computation} \\
C. Jäkel (2023) \cite{jakel2023computation}
\end{tabular} \\
        $\bm{D(10)}$ & \textbf{$\bm{\approx 8.93345 \times 10^{78}}$\ S.E.\ $\bm{2.44 \times 10^{74}}$} & \textbf{This work. See Equation \ref{eq:bestD10Est}} \\
        $\bm{D(11)}$ & \textbf{$\bm{\approx 3.63437\times 10^{144}}$\ S.E.\ $\bm{1.67\times 10^{141}}$} & \textbf{This work. See Table \ref{tbl:bestOtherDededkindEstimations}} \\
        $\bm{D(12)}$ & \textbf{$\bm{\approx 7.14919\times 10^{283}}$\ S.E.\ $\bm{5.35\times 10^{280}}$} & \textbf{This work. See Table \ref{tbl:bestOtherDededkindEstimations}} \\
        $\bm{D(13)}$ & \textbf{$\bm{\approx 6.03589\times 10^{525}}$\ S.E.\ $\bm{1.33\times 10^{523}}$} & \textbf{This work. See Table \ref{tbl:bestOtherDededkindEstimations}} \\
        $\bm{D(14)}$ & \textbf{$\bm{\approx 5.58483\times 10^{1043}}$\ S.E.\ $\bm{4.06\times 10^{1041}}$} & \textbf{This work. See Table \ref{tbl:bestOtherDededkindEstimations}} \\
        $\bm{D(15)}$ & \textbf{$\bm{\approx 3.80603\times 10^{1953}}$\ S.E.\ $\bm{5.30\times 10^{1951}}$} & \textbf{This work. See Table \ref{tbl:bestOtherDededkindEstimations}}
    \end{tabular}
    \caption{Dedekind Numbers from OEIS\cite{oeisA000372} with our estimates for D(10) - D(15) with the smallest Standard Error.}
    \label{tab:dedekindNumbers}
\end{table}

MBFs can be grouped into equivalence classes under permutation of their inputs; exploiting this symmetry is essential for efficient random sampling (Section~\ref{sect:generationLUT7}). Table~\ref{tab:equivalenceClassCounts} lists the known counts of such classes, denoted by $R(n)$.

\begin{table}[H]
    \centering
    \begin{tabular}{c|l|c}
        R(0) & 2 & \\
        R(1) & 3 & \\
        R(2) & 5 & \\
        R(3) & 10 & \\
        R(4) & 30 & \\
        R(5) & 210 & \\
        R(6) & 16353 & \\
        R(7) & 490013148 & T. Stephen \& T. Yusun (2014) \cite{STEPHEN201415} \\
        R(8) & 1392195548889993358 & Bartłomiej Pawelski (2021) \cite{pawelski21} \\
        R(9) & 789204635842035040527740846300252680 & Bartłomiej Pawelski (2023) \cite{pawelski23}
    \end{tabular}
    \caption{Known Equivalence Class Counts, OEIS series A003182 \cite{oeisA003182}}
    \label{tab:equivalenceClassCounts}
\end{table}

In the absence of exact values, asymptotic approximations by Korshunov~\cite{korshunov1981} were the only available estimates for $D(10)$ and beyond. These formulas are highly inaccurate: errors range from 3\% to over 100\% for the known values (Table~\ref{tbl:korshunovEst}), and the formula for odd $n$ is especially unreliable. The consequences are real: when Van Hirtum et al.\ computed $D(9)$, Korshunov's estimate differed from the true value by $-59.72\%$ --- a discrepancy so large that the authors suspected an error in their own code and delayed publication, leading to simultaneous independent publication with~\cite{jakel2023computation}.

This paper provides accurate statistical estimates for $D(10)$ through $D(15)$, with explicit standard errors. We generate large samples of uniformly random MBFs using precomputed lookup tables for $n \leq 7$ and pair matching for larger $n$, and independently via Markov Chain Monte Carlo (MCMC) with degree correction. Three estimation methods are then applied: pair matching frequency, reference subset sampling, and weight layer branching. Estimates from the different methods agree to within their standard errors, providing mutual cross-validation.
All code is available at \url{https://github.com/VonTum/Dedekind} and \url{https://github.com/AlexFihman/Mbf10Up}.

Concurrently with this work, Chen et al.~\cite{chen2026finitenestimatededekindnumbers} independently reported $D(10) = (8.9360 \pm 0.0010) \times 10^{78}$, using a layer-ratio Monte Carlo method related to our weight layer branching approach (Section~\ref{sect:weightLayerBranching}). Their estimate is consistent with our results (Table~\ref{tab:comparisonD10}), providing further independent confirmation.

\section{Korshunov Estimates}
Korshunov\cite{korshunov1981} provides two asymptotic formulas for $D(n)$, one for odd $n$ and one for even $n$. They are not accurate and lack error margins. For odd values the formula is especially inaccurate.

%
%
%
%
%

\begin{table}[H]
    \centering
    \begin{tabular}{c|c|c}
         & Estimate & Error \\\hline
         D(0) & $2.00000$ & $0.00\%$ \\
         D(1) & $6.80423$ & $126.81\%$ \\
         D(2) & $6.59489$ & $9.91\%$ \\
         D(3) & $4.60271\times 10^{1}$ & $130.14\%$ \\
         D(4) & $1.85190\times 10^{2}$ & $10.23\%$ \\
         D(5) & $1.53486\times 10^{4}$ & $102.46\%$ \\
         D(6) & $8.15160\times 10^{6}$ & $4.13\%$ \\
         D(7) & $2.67911\times 10^{12}$ & $10.95\%$ \\
         D(8) & $5.42792\times 10^{22}$ & $-3.30\%$ \\
         D(9) & $1.15330\times 10^{41}$ & $-59.72\%$ \\
         \textbf{D(10)} & $\bm{8.55613\times 10^{78}}$ & $\bm{-4.22\%}$ \\
         \textbf{D(11)} & $\bm{2.56214\times 10^{144}}$ & $\bm{-29.50\%}$ \\
         \textbf{D(12)} & $\bm{6.93452\times 10^{283}}$ & $\bm{-3.00\%}$ \\
         \textbf{D(13)} & $\bm{8.02244\times 10^{525}}$ & $\bm{32.92\%}$ \\
         \textbf{D(14)} & $\bm{5.51832\times 10^{1043}}$ & $\bm{-1.19\%}$ \\
         \textbf{D(15)} & $\bm{7.93171\times 10^{1953}}$ & $\bm{108.40\%}$ \\
    \end{tabular}
    \caption{Korshunov estimates for Dedekind numbers with relative error to the values from Table \ref{tab:dedekindNumbers}. For D(0)--D(9) the error is relative to exact known values; for D(10)--D(15) (bold) it is relative to this paper's estimates (See Table \ref{tab:dedekindNumbers}) and is therefore itself uncertain.}
    \label{tbl:korshunovEst}
\end{table}


\section{Methods for Generating Uniformly Random Monotone Boolean Functions}

We describe methods for generating uniformly random MBFs: lookup tables for small dimensions ($n \leq 7$), and recursive pair matching for higher dimensions. The resulting samples serve both as a source for direct estimation and as a reference for validating the MCMC method.

\subsection{Lookup Tables for \texorpdfstring{$n \leq 6$}{n <= 6}}
\label{sect:generationLUT6}

A Boolean Function of $n$ variables can be represented as a $2^n$ bit bitset in computer memory. Each bit in the bitset represents the output of the function when the inputs of the function match the binary representation of the bit's index, like a boolean truth table. 
With this encoding scheme, the largest of these lookup tables that can be stored on contemporary hardware is for MBFs in 6 variables, also referred to as MBF6s. A list of all distinct 64-bit MBF6s ($D(6)=7828354$ entries) takes up 63MB.
Random sampling from this precomputed list yields a stream of uniformly sampled MBF6s. 

Doing the same for MBF7 would require $D(7)=2414682040998$ 128-bit elements, or 36TB of memory. This falls outside the capacity of any system at our disposal. 

\subsection{Lookup Table for \texorpdfstring{$n = 7$}{n = 7}}
\label{sect:generationLUT7}
As will be mentioned in Section \ref{sect:generationPairMatching}, a lot of unneeded MBF6 samplings could be saved if we could directly sample MBF7s instead of combining MBF6. 

Since the full MBF7 lookup table would be too large, we need to store it in a more efficient encoding. We can do this by exploiting the symmetries inherent in MBFs. MBFs can be categorized into equivalence classes over permutation of the function inputs. If MBF $\alpha$ can be turned into MBF $\beta$ by permuting its inputs, (such as swapping input 1 and 3), they belong to the same equivalence class and are therefore stored under the same "representative". For $n=7$ the number of equivalence classes is $R(7)=490013148$ (See Table \ref{tab:equivalenceClassCounts}). If we store only one representative per equivalence class, we reduce our memory footprint to roughly 8GB. The sizes of these equivalence classes are divisors of $n!$, with the vast majority of them being equal to $n!$. A dataset of all inequivalent MBF7s was computed in \cite{vanhirtum2023computation}, which made it available at \url{https://zenodo.org/records/10579803}. 

Because the equivalence classes are not all of the same size, simple naive sampling into this buffer will skew our distribution towards MBFs with high symmetry. The MBFs from smaller classes would be overrepresented in the final distribution. To compensate for this skew, we will weight the probability of selecting each equivalence class by its size. Practically, we do this by sorting the equivalence classes by size first. The groupings of different Equivalence Class size for 7 variables are shown in Table \ref{tab:mbfsByClassSize}. 

\begin{table}[H]
    \centering
    \begin{tabular}{r|l c r|l}
        Class Size & Classes of this size && Class Size & Classes of this size \\
        \cmidrule(lr){1-2}\cmidrule(lr){4-5}
        5040 & 468822749 && 210 & 3255 \\
        2520 & 20005503  && 140 & 702 \\
        1680 & 3128      && 120 & 4 \\
        1260 & 1024050   && 105 & 1206 \\
        840  & 75024     && 84  & 9 \\
        720  & 4         && 70  & 90 \\
        630  & 47242     && 42  & 99 \\
        504  & 237       && 35  & 117 \\
        420  & 26739     && 30  & 5 \\
        360  & 18        && 21  & 75 \\
        315  & 2742      && 7   & 27 \\
        252  & 114       && 1   & 9 \\
    \end{tabular}
    \caption{Number of equivalence classes by Equivalence Class Size for MBFs of 7 Variables. The second column sums to $R(7)$. The sum of the products of the columns gives $D(7)$. See \cite{vanhirtum2023computation}.}
    \label{tab:mbfsByClassSize}
\end{table}

To use this table, the sampling algorithm selects a random number $r$ between $0$ and $D(n)$. Then going through the table, for each entry it checks if $r < ClassSize \cdot \#ClassesOfSize$. If not, it subtracts $ClassSize \cdot \#ClassesOfSize$ from $r$ and moves to the next row, until it finds one that is smaller. Because the first entry of $ClassSize = n!$ accounts for the vast majority of the probability space, it gets chosen the bulk of the time, so having it checked first improves performance. 

Once we have selected a buffer, we floor-divide $r / ClassSize$
to derive the MBF index in this buffer. Finally, permuting the MBF randomly results in a uniformly sampled MBF7. 

Further optimizations to this algorithm were made, such as reducing soft page faults by enabling Huge Pages, and prefetching the next 64 MBFs to pipeline the DRAM access latency. 


\subsection{Generating MBFs through Pair Matching}
\label{sect:generationPairMatching}

The following Lemma (Lemma \ref{lemma:pairMatch}) is used to generate a uniform sampling of $n + 1$ variable MBFs, given a uniform sampling stream of $n$ variable MBFs. What it describes is that given two uniformly randomly sampled MBFs in $n$ variables, they can be checked for a "match" against Equation \ref{eq:matchCriterion}. If they adhere to it, they can be combined to form a new (now uniformly randomly sampled) MBF in $n+1$ variables. If not, both MBFs sampled must be discarded. Reusing previously generated MBFs skews the output distribution in favor of MBFs with lower match probability. 

\begin{lemma}
Let $f_1, f_2 : \{0,1\}^n \to \mathbb{R}$ be two boolean functions such that
\begin{equation}
\forall \tilde{\alpha} \in \{0,1\}^n,\quad f_1(\tilde{\alpha}) \leq f_2(\tilde{\alpha}).
\label{eq:matchCriterion}
\end{equation}
Then the pair $(f_1,f_2)$ defines a function $f \in M(n+1)$ by
\begin{equation}
f(x_0, x_1, \dots, x_n) =
\begin{cases}
f_1(x_1, \dots, x_n), & \text{if } x_0 = 0, \\
f_2(x_1, \dots, x_n), & \text{if } x_0 = 1.
\end{cases}
\label{eq:matchCombine}
\end{equation}
\label{lemma:pairMatch}
\end{lemma}
\begin{corollary}
As a corollary of Lemma \ref{lemma:pairMatch}, it follows that any MBF in $n$ variables can be split into a unique pair of MBFs in $n - 1$ variables, for which Equation \ref{eq:matchCriterion} holds. 
\label{corollary:splitIntoPairs}
\end{corollary}

The main bottleneck of pair matching is the large number of lower-order MBFs required to generate a single higher-order MBF. On average 50.7 MBF6 are required to generate an MBF7, 207.8 MBF7 for an MBF8, and 22002.6 MBF8 for an MBF9. Multiplying these together, generating a single MBF9 would require on average 232 million MBF6. By using the MBF7 lookup table from Section~\ref{sect:generationLUT7}, the number of samples required per MBF9 is reduced by a factor of 50 to only 4.6 million MBF7.

The original MBF6 based implementation averaged 0.4 MBF9 per second per thread. The MBF7 based implementation generated at 860 MBF9 per second on a 128 core Noctua 2 Compute Node \cite{noctua2}. Per-thread comparison is difficult since the one is L3 cache bound, and the other is memory bandwidth bound. We ran this to generate a $2.3 \cdot 10^9$ MBF9 dataset which we reuse in later sections. 

\subsection{Generating MBFs by Random Walks}
\label{sect:walks}

We can look at the body of Monotone Boolean functions of certain dimension as an undirected graph, $G(n) = G(M(n), E(n))$, where nodes represent functions in $M(n)$, and edges $E(n)$ connect functions that differ only in one input value. In other words, edges connect functions whose Hamming distance is exactly 1.

The graph $G(n)$ is \textit{connected}: every non-zero MBF has at least one value that can be changed from 1 to 0 without violating monotonicity. Repeating this process eventually reaches the constant-0 function. Similarly, every non-one MBF can reach the constant-1 function by repeatedly changing suitable values from 0 to 1.

For any $n$, the graph $G(n)$ is \textit{bipartite}. It can be split into two parts: nodes of odd and even weight. Odd-weight nodes are connected to even-weight nodes only, and vice versa.

\subsubsection{Finding nearby nodes}
Consider some $f \in M(n)$. We define the \textbf{Minimal Positive Set} of $f$ as the set of minimal elements where $f$ has value 1:

\begin{equation}
\operatorname{MinPos}(f) = \{x \in \{0,1\}^n : f(x) = 1 \text{ and } \forall y < x, f(y) = 0\}
\end{equation}

Similarly, we define the \textbf{Maximal Negative Set} as the maximal elements where $f$ has value 0:

\begin{equation}
\operatorname{MaxNeg}(f) = \{x \in \{0,1\}^n : f(x) = 0 \text{ and } \forall z > x, f(z) = 1\}
\end{equation}
where $y<x$ denotes the coordinate-wise partial order on the Boolean hypercube.

These two sets completely characterize the monotone Boolean function, since by monotonicity:
\begin{itemize}
\item $\forall x \in \operatorname{MinPos}(f), \forall z \geq x: f(z) = 1$
\item $\forall x \in \operatorname{MaxNeg}(f), \forall y \leq x: f(y) = 0$
\end{itemize}

The sets $\operatorname{MinPos}(f)$ and $\operatorname{MaxNeg}(f)$ provide an alternative representation of an MBF and correspond to its minimal DNF and CNF representations, respectively. Both sets are antichains: if a point belongs to either set, no comparable point can belong to the same set. As a result, they tend to occupy the middle layers of the Boolean hypercube, which contain the largest number of mutually incomparable vertices.

For even dimension $n$, the middle layer at weight $n/2$ contains the most vertices. For odd dimension $n$, the two middle layers at weights $(n-1)/2$ and $(n+1)/2$ do so. Since these layers can accommodate the largest antichains, in uniformly sampled monotone Boolean functions, most elements of $\operatorname{MinPos}(f)$ and $\operatorname{MaxNeg}(f)$ are found in these layers.

As the dimension grows, this effect becomes stronger. Up to $n=9$, the minimal sets are distributed across several layers, although layers 4 and 5 already dominate at $n=9$. For $n>9$, a large majority of the elements of $\operatorname{MinPos}(f)$ and $\operatorname{MaxNeg}(f)$ are concentrated in a single middle layer, while the layers below and above become nearly saturated with 0s and 1s, respectively. This creates two symmetric regions in the graph $G(n)$, one for each choice of middle layer. By duality, both regions contain the same number of MBFs, so no special treatment is required in the computations. However, the connectivity between these regions becomes increasingly poor as the dimension grows.

For an MBF, the value at a single input can be flipped without violating monotonicity only if that input belongs to $\operatorname{MinPos}(f)$ or $\operatorname{MaxNeg}(f)$. These two sets therefore determine all function's neighbors. In the MCMC algorithm, we maintain them explicitly and update them after each step.

\subsubsection{Random Walk and Degree Correction}
\label{sect:degreeCorrection}
Starting from some arbitrary MBF, flipping a random bit in $\operatorname{MinPos}(f)$ or $\operatorname{MaxNeg}(f)$ produces a new neighboring MBF. We define $\deg(f) = |\operatorname{MinPos}(f)| + |\operatorname{MaxNeg}(f)|$ as the number of neighbors of $f$ in the MBF graph. Because the graph is connected, finite and bipartite, there are stationary distributions for odd and for even nodes. For a simple random walk on an undirected graph, the stationary distribution is $\pi(f) \propto \deg(f)$, so high-degree nodes are visited more frequently than low-degree ones. This means that raw MCMC samples are not uniformly distributed over $M(n)$; whenever a uniform average is required, each sample $f$ must be weighted by $1/\deg(f)$ to correct for this bias.

Since consecutive MCMC samples are correlated, they cannot be treated as independent draws. Nevertheless, we use MCMC to generate random MBFs throughout this work. We focus on dimensions 10 to 15, while dimension 9 is used to validate the method against uniform samples. Within each densely populated region, the graph is well connected, since typical MBFs have many flippable inputs in the middle layers. The main bottleneck is between the two symmetric regions that arise in higher odd dimensions. 

For any two MBFs $f,g \in M(n)$, the graph distance between them is exactly their Hamming distance. Indeed, consider an input $x$ where $f(x)\neq g(x)$. If $f(x)=1$ and $g(x)=0$, then among the remaining disagreements of this type there is a minimal one; flipping it from 1 to 0 preserves monotonicity. Similarly, if $f(x)=0$ and $g(x)=1$, then among the remaining disagreements of this type there is a maximal one; flipping it from 0 to 1 preserves monotonicity. Repeating this process changes one disagreeing input at a time and reaches $g$ after exactly $d_H(f,g)$ steps. Thus, the shortest path between two MBFs in $G(n)$ has length equal to their Hamming distance, which is at most $2^n$.

The convergence to the stationary distribution occurs at an exponential rate. After $N$ steps (known as the mixing time),  the total variation distance from the stationary distribution decreases by the factor of $\epsilon = \frac{1}{4}$. (Theorem 4.3 \cite[p.~52]{levin2017markov})

For practical purposes, we estimate the number of steps after which the deviation from the stationary distribution falls below a detectable threshold. We used two indirect methods to assess this.

First, we compare the average Hamming distance between two uniformly random MBFs and MBFs obtained through MCMC (Section \ref{sect:hamming9}). Second, we analyze the weight distribution. For uniformly random MBFs, the weight distribution is symmetric about the midpoint $2^{n-1}$ (Section \ref{sect:weightBalance}).

\subsubsubsection{Expected Hamming distance, \texorpdfstring{$M(9)$}{M(9)}}
\label{sect:hamming9}

The expected \textit{Hamming distance} between two uniformly random MBFs was calculated from the MBF9s generated using the method in Section \ref{sect:generationPairMatching}. We sampled $167,772,160$ random MBF9. Result: $mean = 86.151734$, $\sigma = 11.432587$, $S.E. = 8.8 \times 10^{-4}$.

The starting point for our MCMC chains is the function 0. In order to compensate for odd/even cycle, we take $k$ or $k+1$ steps with equal probability. Table \ref{tab:average_hamming_distance9} shows the statistics of the Hamming distances that we found. After about 15,000 iterations, we see convergence of $mean$ and $\sigma$. This is illustrated in Table \ref{tab:average_hamming_distance9}. 

\begin{table}[H]
    \centering
    \begin{tabular}{c|c|c}
        \toprule
        \textbf{Steps} & \textbf{Mean distance} & $\boldsymbol{\sigma}$ \\
        \midrule
        100   & 31.8427903  & 5.868380023  \\
        300   & 52.49736939 & 7.020338912  \\
        1000  & 74.92570711 & 7.923565358  \\
        3000  & 81.36795606 & 9.288982428  \\
        5000  & 84.68916137 & 10.7902874   \\
        10000 & 86.0865943  & 11.40300854  \\
        15000 & 86.14450992 & 11.42455114  \\
        20000 & 86.15614547 & 11.43540288  \\
        30000 & 86.15729709 & 11.43944073  \\
        60000 & 86.13405866 & 11.43576622  \\
        \bottomrule
    \end{tabular}
    \caption{Average Hamming distance by number of steps for $M(9)$. The experiment was repeated $10^6$ times.}
    \label{tab:average_hamming_distance9}
\end{table}

\subsubsubsection{Hamming distance, \texorpdfstring{$M(10)$}{M(10)}}
\label{sect:hamming10}
The same convergence analysis was performed on the $M(10)$ graph. Although we do not have a source of uniformly random MBFs of this dimension, we can compare the $mean$ and $\sigma$ against a long-running MCMC chain.

\begin{table}[H]
    \centering
    \begin{tabular}{c|c|c}
        \toprule
        \textbf{Steps} & \textbf{Mean distance} & $\boldsymbol{\sigma}$ \\
        \midrule
        100   & 37.77584935 &  6.566062941 \\
        300   & 69.0808887  &  8.339321822 \\
        1000  & 113.982995  &  10.07327182 \\
        3000  & 141.5978466 &  12.52099822 \\
        5000  & 151.9950193 &  14.30656889 \\
        10000 & 141.609136  &  9.417690752 \\
        15000 & 141.5688093 &  9.394752421 \\
        20000 & 141.55493   &  9.37813645  \\
        30000 & 141.5590908 &  9.390921352 \\    
        60000 & 141.5461016 &  9.38908007  \\
        \bottomrule
    \end{tabular}
    \caption{Average Hamming distance by number of steps for $M(10)$. The experiment was repeated $10^6$ times.}
    \label{tab:average_hamming_distance10}
\end{table}

Here we see convergence of the $mean$ and $\sigma$ after about 15,000 iterations.
A long-running experiment of $10^9$ steps on two independent MCMC chains gives $mean = 141.553$ and $\sigma = 9.3839$.

\subsubsubsection{MBF Weight balance}
\label{sect:weightBalance}
Another measure of convergence is the balance of weight distribution.

For a function in $M(n)$, there are $2^n + 1$ possible weight values. The weight distribution of MBFs in $M(n)$ is by definition symmetrical (due to each MBF of weight $w$ corresponding to a unique dual MBF of weight $2^n - w$), so in a uniform sample the number of MBFs of weight $w$ should roughly equal the number of functions of weight $2^n-w$. An example of a poor weight distribution (5{,}000 steps) alongside a converged one is shown in Figure \ref{fig:mbf9DistributionComparison}.

To compute the Imbalance, we define $L$ and $R$ as the sums of the degree-corrected weight counts $P_i$ (Equation \ref{eq:Pi}) below and above the midpoint $2^{n-1}$:

\begin{equation}
    \label{eq:leftSum}
    L = \sum_{i=0}^{2^{n-1}-1} P_i
\end{equation}

\begin{equation}
    \label{eq:rightSum}
    R = \sum_{i=2^{n-1}+1}^{2^{n}} P_i
\end{equation}

\begin{equation}
    \label{eq:imbalance}
    Imbalance = \frac{\left| L - R \right|}{L + R}
\end{equation}

where $P_i$ is the degree-corrected count of sampled MBFs with weight $i$: (see Section~\ref{sect:walks} for the correction rationale)
\begin{equation}
    \label{eq:Pi}
    P_i = \sum_{\substack{f \in \text{sample} \\ |f| = i}} \frac{1}{\deg(f)},
\end{equation}

\begin{table}[H]
    \centering
    \begin{tabular}{c|c}
        \toprule
        \textbf{Steps} & \textbf{Imbalance} \\
        \midrule
        5000 & $3.98 \times 10^{-1}$\\
        10000 & $7.13 \times 10^{-2}$ \\
        15000 & $1.33 \times 10^{-2}$ \\
        20000 & $2.63 \times 10^{-3}$ \\
        25000 & $7.68 \times 10^{-4}$ \\
        30000 & $3.95 \times 10^{-4}$ \\
        \bottomrule
    \end{tabular}
    \caption{LR-Imbalance by number of steps, on the $M(9)$ graph}
    \label{tab:lr_imbalance9}
\end{table}

For comparison, the expected imbalance due to pure sampling fluctuation in a perfectly balanced distribution of $10^7$ samples is $2.52 \times 10^{-4}$, derived as follows. Let $L$ and $R$ count samples falling on each side; the imbalance is:  

\begin{equation}
\text{Imbalance} = \frac{|L - R|}{L + R}
\end{equation}

Since \( L \sim \text{Binomial}(10^7, 0.5) \), using the normal approximation:

\begin{equation}
L - R \sim \mathcal{N}(0, 10^7), \quad \sigma = \sqrt{10^7} \approx 3162.28
\end{equation}

The expected absolute difference:

\begin{equation}
\mathbb{E}[|L - R|] \approx 3162.28 \times \sqrt{\frac{2}{\pi}} \approx 2523.1
\end{equation}

Thus, the expected imbalance is:

\begin{equation}
\mathbb{E}[\text{Imbalance}] \approx \frac{2523.1}{10^7} \approx 0.000252 
\end{equation}

The same convergence analysis was performed on the $M(10)$ graph, where convergence is faster.
\begin{table}[H]
    \centering
    \begin{tabular}{c|c}
        \toprule
        \textbf{Steps} & \textbf{Imbalance} \\
        \midrule
        5000 & $8.98 \times 10^{-1}$\\
        10000 & $2.48 \times 10^{-2}$ \\
        15000 & $9.2 \times 10^{-5}$ \\
        \bottomrule
    \end{tabular}
    \caption{LR-Imbalance by number of steps, on the $M(10)$ graph}
    \label{tab:lr_imbalance10}
\end{table}

\subsection{Validation of MBF9 generated by different methods}

\begin{figure}[H]
    \centering
    \includegraphics[width=1.0\linewidth]{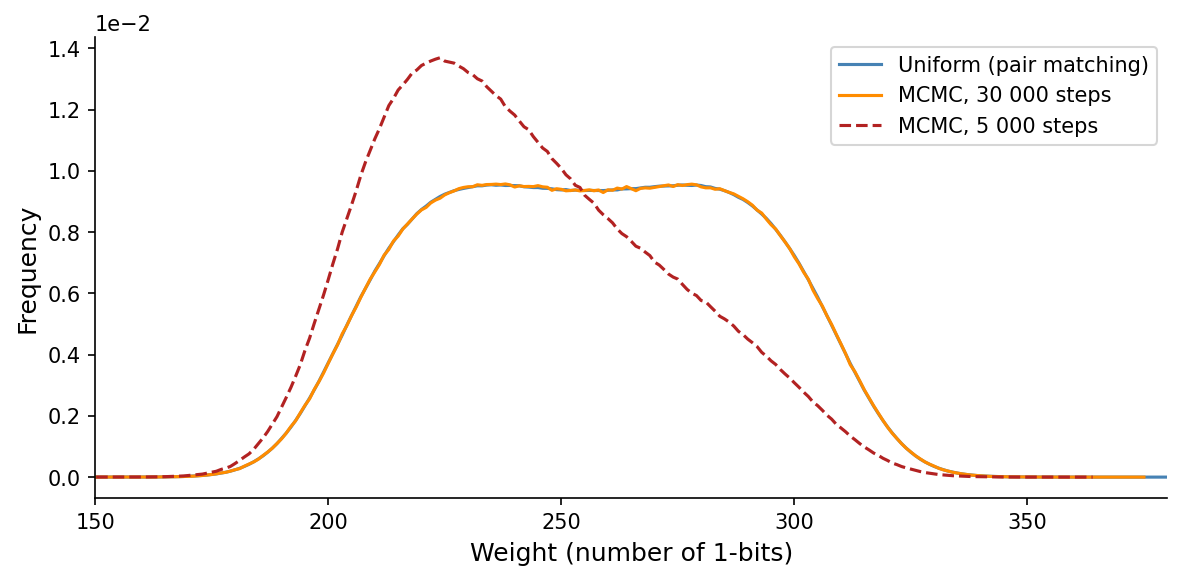}
    \caption{Degree-corrected weight distribution of MBF9s from MCMC at 5{,}000 and 30{,}000 burn-in steps, compared against the uniform pair-matching reference ($1.56 \times 10^8$ samples). The bimodal shape of the uniform distribution reflects the two middle layers (4 and 5) of the Boolean hypercube. After 30{,}000 steps the MCMC distribution is statistically indistinguishable from the uniform reference (KS $p = 0.29$).}
    \label{fig:mbf9DistributionComparison}
\end{figure}

Since the MCMC stationary distribution is $\pi(f) \propto \deg(f)$, the degree-corrected weight histogram (each sample weighted by $1/\deg(f)$; see Section~\ref{sect:walks}) should match the uniform pair-matching distribution if the chain has converged. We test this directly with the two-sample Kolmogorov--Smirnov test on the weight histograms. With $10^7$ chains (30{,}000 burn-in steps) against $1.56 \times 10^8$ pair-matched samples, the KS statistic is $D = 0.000320$ ($p = 0.29$): no evidence of a distributional difference. After only 5{,}000 steps, $D = 0.197$ ($p \approx 0$), confirming that the shorter burn-in is insufficient. Figure~\ref{fig:mbf9DistributionComparison} shows both MCMC distributions alongside the uniform reference.

\section{Methods for estimating Higher Dedekind Numbers}

\subsection{Pair Matching}

Theorem \ref{theorem:pairMatchEstimation} links $D(n)$ and $D(n+1)$ through the pair-matching probability.

\begin{theorem}
Due to Corollary \ref{corollary:splitIntoPairs} any MBF in $n + 1$ variables has exactly one split into a unique pair of MBFs in $n$ variables. This means that the probability that a uniformly sampled pair of MBFs in $n$ variables will successfully combine is:
\begin{equation}
    p_{match, n} = \frac{D(n+1)}{D(n)^2}
    \label{eq:computePMatchFromDedek}
\end{equation}

Therefore given accurate estimates for $D(n)$ and $p_{match, n}$, $D(n+1)$ can be estimated as:
\begin{equation}
D(n+1) = D(n)^2 \cdot p_{match, n}
    \label{eq:computeDedekFromPMatch}
\end{equation}

The best currently known values of $p_{match, n}$ are provided in Table \ref{tab:combinationProbabilities}.
\begin{table}[H]
    \centering
    \begin{tabular}{c|l|l}
        $p_{match,n}$ & Value & S.E. \\\hline
        $p_{match,0}$ & 3/4 & 0 \\
        $p_{match,1}$ & 2/3 & 0 \\
        $p_{match,2}$ & 5/9 & 0 \\
        $p_{match,3}$ & 0.42 & 0 \\
        $p_{match,4}$ & 0.26860119047619047 & 0 \\
        $p_{match,5}$ & 0.13621265655199447 & 0 \\
        $p_{match,6}$ & 0.03940207011036327 & 0 \\
        $p_{match,7}$ & 0.00962672400664057 & 0 \\
        $p_{match,8}$ & 0.00009089831010396 & 0 \\
        $p_{match,9}$ &  0.00010892156976784  &  $2.97\cdot 10^{-09}$ \\
        $p_{match,10}$ &  $4.55397\cdot 10^{-14}$  &  $2.10\cdot 10^{-17}$ \\
        $p_{match,11}$ &  $5.41250\cdot 10^{-06}$  &  $6.41\cdot 10^{-09}$ \\
        $p_{match,12}$ &  $1.18093\cdot 10^{-42}$  &  $3.14\cdot 10^{-45}$ \\
        $p_{match,13}$ &  $1.53294\cdot 10^{-08}$  &  $1.30\cdot 10^{-10}$ \\
        $p_{match,14}$ &  $1.22025\cdot 10^{-134}$  &  $2.45\cdot 10^{-136}$ \\
    \end{tabular}
    \caption{Best currently known values of $p_{match,n}$. The values up to $n \leq 8$ are known exactly from known Dedekind Numbers (See Equation \ref{eq:computePMatchFromDedek}), higher values are derived backwards from our best $D(n)$ estimates from Table \ref{tab:dedekindNumbers}. }
    \label{tab:combinationProbabilities}
\end{table}

\label{theorem:pairMatchEstimation}
\end{theorem}

Now that we have a stream of randomly sampled MBF9s, we can try to combine them to create MBF10s as we did in Section \ref{sect:generationPairMatching} and record the number of successes and failures. Performing this for a large number of such attempts allows us to accurately estimate the ratio of successful comparisons. In the case for combining two MBF9s to create MBF10s, a successful combination occurs about once every 10000 comparisons. (See Equation \ref{eq:firstEstimatedMatchProbability}). 

This probability is sufficient for practical computation, and yields an estimate better than Korshunov's asymptotic formula. Each comparison is a Bernoulli trial (resulting in either success or failure), making the sample mean an unbiased Monte Carlo estimator of the underlying probability \cite{jiang2014guaranteedmontecarlomethods}.

Given this, about $10000 \left( \frac{1}{p_{match, 9}} - 1 \right) \approx 10^8$ comparisons are required to achieve a relative Standard Error of $10^{-2}$.

\subsection{Estimating \texorpdfstring{$D(10)$}{D(10)} using MBF9 Buffer matching}
\label{sect:leftRightMBFDefinition}

Generating uniformly random MBF9s (through the pair matching method from Section \ref{sect:generationPairMatching}) is very expensive. To maximize the reuse of MBF9s that we generate, we create two large buffers of MBF9s called \texttt{left} and \texttt{right}. By doing this we can perform $|left||right|$ matches, instead of $(|left|+|right|) / 2$. 

In a preliminary experiment we performed $2.5 \times 10^7$ loops, using $10^4$ \texttt{left} MBFs and $10^4$ \texttt{right} MBFs, resulting in the estimation:
\begin{equation}
\hat{p}_{match, 9} = 1.08808072 \times 10^{-4}, S.E. = 2.8 \times 10^{-7}
\label{eq:firstEstimatedMatchProbability}
\end{equation}

which using Equation \ref{eq:computeDedekFromPMatch} results in $D(10) = 8.92(41) \times 10^{78}$ with $\sigma = 2.3 \times 10^{76}$.

For comparison, without buffering the same $5 \times 10^7$ raw MBF9s would be used as $2.5 \times 10^7$ non-overlapping pairs. Each pair is a Bernoulli trial, so:
\[
\sigma_{\text{non-buf}}(D(10)) = D(9)^2 \sqrt{\frac{p}{N}} = D(9)^2 \sqrt{\frac{1.08808 \times 10^{-4}}{2.5 \times 10^7}} \approx 1.71 \times 10^{77}
\]

By reusing the same MBF9s through buffering, $\sigma$ is reduced by a factor of $7.47$ to $2.3 \times 10^{76}$, corresponding to a speedup of $55.9\times$.

The total match count from the $|L| \times |R|$ comparison matrix is a sum of dependent Bernoulli variables, but the CLT still applies. To see this, note that each cyclic super-diagonal of the matrix — the set of pairs $\{(i,\, i+k \bmod |R|)\}_{i=0}^{|L|-1}$ for a fixed offset $k$ — consists entirely of independent comparisons (distinct left and right MBFs). The total count is a sum of such diagonals, each binomially distributed, so the overall distribution has bounded higher-order moments and well-behaved variance.

To characterize the statistical properties of $\hat{p}_{match,9}$ empirically, we performed two block-comparison experiments on $1.56 \times 10^8$ uniformly random MBF9s: blocks of $10^5 \times 10^5$ (838 trials) and $10^4 \times 10^4$ (8388 trials). Both give a consistent mean $\bar{p} = 1.088 \times 10^{-4}$. The distribution of $\hat{p}$ is approximately Gaussian --- at the $10^4$ scale we measure skewness $0.28$ and excess kurtosis $0.22$, which collapse toward zero at the $10^5$ scale by the Central Limit Theorem. The observed standard deviation substantially exceeds the Bernoulli prediction: $\sigma / \sigma_{\text{Bernoulli}} = 9.47$ and $30.5$ respectively, reflecting that comparisons within a block are correlated. The effective number of independent comparisons scales as $N_{\text{eff}} \propto B$ rather than $B^2$, where $B$ is the block size --- each additional MBF contributes a constant amount of independent information regardless of how many cross-comparisons it participates in. The mild excess kurtosis ($0.22 \pm 0.05$) has negligible effect on the mean estimate: the correction to the standard error of the mean is less than $0.001\%$.

The exhaustive $|L| \times |R|$ comparison has $O(B^2)$ complexity in the block size $B$; implementation optimizations --- a bitset transpose, AVX2 SIMD matching, and a binary-tree pre-filter that eliminates 92\% of comparisons at the largest block size --- are described in Appendix~\ref{app:bufferImpl}.

\subsubsection{Results from Pair Matching with Filter Tree}

\begin{table}[H]
\centering
\begin{tabular}{c|c|c|c|c|c}
Block Size & Trials & {Mean Match} & {$D(10)$} & {Std.\ Dev.} & {S.E.} \\
           &         & {Fraction} & {Estimate} & {(D10)} & {(Calculated)} \\
\midrule
100000    & 23655 & 0.000108929 & 8.93402e+78 & 3.60967e+77 & 2.347e+75 \\
300000    & 7885 & 0.000108944 & 8.93532e+78 & 2.08813e+77 & 2.352e+75 \\
1000000   & 2365 & 0.000108957 & 8.93636e+78 & 1.16641e+77 & 2.398e+75 \\
3000000   & 788 & 0.000108955 & 8.93618e+78 & 6.46594e+76 & 2.303e+75 \\
10000000  & 236 & 0.000108975 & 8.93782e+78 & 3.4312e+76  & 2.234e+75 \\
\end{tabular}

\caption{Estimates for $D(10)$ using pair matching on various block sizes, using the 2.3B MBF9 dataset generated in Section \ref{sect:generationPairMatching}. The block size refers to the total number of MBFs involved in a single buffer comparison trial. As an example, if the block size is 100k, then that means we compare 50k left MBFs with each of 50k right MBFs, for a total of 2.5B comparisons. Using our 2.3B MBF dataset, we repeat this 23655 times. Compute time (single-core) for pair matching spanned from 88s for Block Size 100k, to 1392s for Block Size 10M. }
\label{tab:pairMatchingByBlockSize}
\end{table}

The gains in efficiency from larger block sizes were limited by the high throughput of the SIMD bitset comparison and by the diminishing per-MBF information as block size increases.

The buffer-based pair matching approach achieves a best S.E. of $2.234 \times 10^{75}$ (block size 10M), which is within one order of magnitude of the best method (S.E. $= 2.44 \times 10^{74}$, Section~\ref{sect:wellKnownMBF}). It therefore provides a meaningful independent estimate and a useful cross-check on the other methods. 

\subsection{Estimating \texorpdfstring{$D(10)$--$D(12)$}{D(10)--D(12)} via Reference Subset Sampling}
\label{sect:wellKnownMBF}

A \emph{reference subset} $S \subseteq M(n)$ is a subset whose cardinality $|S|$ is given by a simple closed-form expression. Given such a subset, the fraction $p_S = |S|/D(n)$ can be estimated by measuring how often MCMC samples fall into $S$, yielding
\[
D(n) = \frac{|S|}{\hat{p}_S}.
\]

Because consecutive MCMC samples are correlated, the standard error of $\hat{p}_S$ cannot be estimated from a single experiment using the Bernoulli formula. Instead, we repeat the experiment $K$ independent times, each generating $N$ MCMC samples and producing one estimate $\hat{p}_S^{(i)}$. The standard error is then:
\[
\text{S.E.}(\bar{p}_S) = \frac{\mathrm{std}(\hat{p}_S^{(i)})}{\sqrt{K}}.
\]

The method is practical when the expected number of hits per experiment $N \cdot p_S \gg 1$, so that each experiment yields a meaningful estimate.

For $n=10$, a natural choice is the set of \emph{1-layer MBFs} at the middle layer. A \textbf{1-layer MBF} at layer $k$ is a monotone Boolean function satisfying $f(x)=0$ for all $|x|<k$ and $f(x)=1$ for all $|x|>k$, with $f$ unrestricted on inputs of weight $k$. Since all $\binom{n}{k}$ inputs at weight $k$ are mutually incomparable under the coordinate-wise partial order, and thus can be added or removed independenly. This means there are exactly $2^{\binom{n}{k}}$ such functions. For $n=10$, $k=5$:
\[
|S| = 2^{\binom{10}{5}} = 2^{252} \approx 7.237 \times 10^{75}.
\]

Using the MCMC algorithm described in Section~\ref{sect:walks}, we generate a sequence of MBF10s. We skip the first 30,000 functions as burn-in, and subsequently record how many belong to this 1-layer subset (layer 5).

In a single experiment we generate $10^8$ functions, and find the sample proportion $\hat{p}_{L5}$ of a function to be of layer 5 only. A degree correction is applied as described in Section~\ref{sect:walks}. From here, 

\begin{equation}
D(10) = \frac{2^{\binom{10}{5}}}{p_{L5}}
\end{equation}

where $p_{L5}$ is the probability that an MBF10 is a 1-layer function at layer 5.

Each experiment takes about 6.5 seconds CPU time, single thread, and 40,000 experiments were performed, giving $\widehat{1/p_{L5}} = 1234.414$ and
\begin{equation}
D(10) \approx 8.93345 \times 10^{78}, S.E. = 2.44 \times 10^{74}
\label{eq:bestD10Est}
\end{equation}

For odd $n$, there are two middle layers. For $n=11$ the middle layers are at weights 5 and 6, each contributing $2^{\binom{11}{5}} = 2^{462}$ 1-layer MBFs. By the duality symmetry of $M(n)$ the two counts are equal, so $|S| = 2^{463}$ and the same estimator applies with $p_S$ now being the fraction of MBF11s that are 1-layer at either middle layer. For $n=12$ (even), the single middle layer at weight 6 gives $|S| = 2^{\binom{12}{6}} = 2^{924}$.

Similar experiments were performed for dimensions 11 and 12, giving the following results:
\begin{equation}
D(11) \approx 3.6325 \times 10^{144}, S.E. = 7.16 \times 10^{141}
\end{equation}
\begin{equation}
D(12) \approx 7.1911 \times 10^{283}, S.E. = 3.25 \times 10^{281}
\end{equation}

For $n \geq 13$ this method becomes impractical. At $n=13$ the estimated $p_S \approx 1.7 \times 10^{-9}$, yielding fewer than one hit per $10^8$-sample experiment on average. The method is therefore not applied beyond $n=12$.

\subsection{Estimating \texorpdfstring{$D(10)$--$D(15)$}{D(10)--D(15)} by Weight Layer Branching}
\label{sect:weightLayerBranching}

In the following experiments, a set of monotone Boolean functions of a certain dimension is considered as a layered graph by weight layers.
\begin{definition}[Weight layer]
A \emph{weight layer} is the subset of the set of monotone Boolean functions of dimension $n$ having a given weight~$x$. That is, in exactly $x$ of the $n^2$ possible inputs, the function returns $1$. 
\end{definition}

In a graph of MBFs, edges exist only between vertices in consecutive layers, i.e., between layers $x$ and $x\pm 1$. An experiment consists of two independent sequences of random walks between nearby weight layers: one beginning at the unique vertex in layer~$0$, and one beginning at the unique vertex in the final layer of weight~$2^n$. For each visit to a vertex in layer~$x$, we record the number of its neighbors in layer~$x+1$ (right degree) and in layer~$x-1$ (left degree). Averaging these values over many visits, with correction for the probability of visiting each node, we get an estimate of the typical branching factor between adjacent layers. The standard degree correction (Section~\ref{sect:degreeCorrection}) is applied: each observation at node $v$ is weighted by $1/\deg(v)$, where $\deg(v)$ is the degree within the local sub-graph of the two adjacent weight layers being traversed — the sum of the left degree (edges to layer $x-1$) and the right degree (edges to layer $x+1$).

\subsubsection{Computation of cardinality of layers and midpoint discrepancy}

\begin{figure}
    \centering
    \begin{subfigure}{0.48\linewidth}
        \centering
        \includegraphics[width=\linewidth]{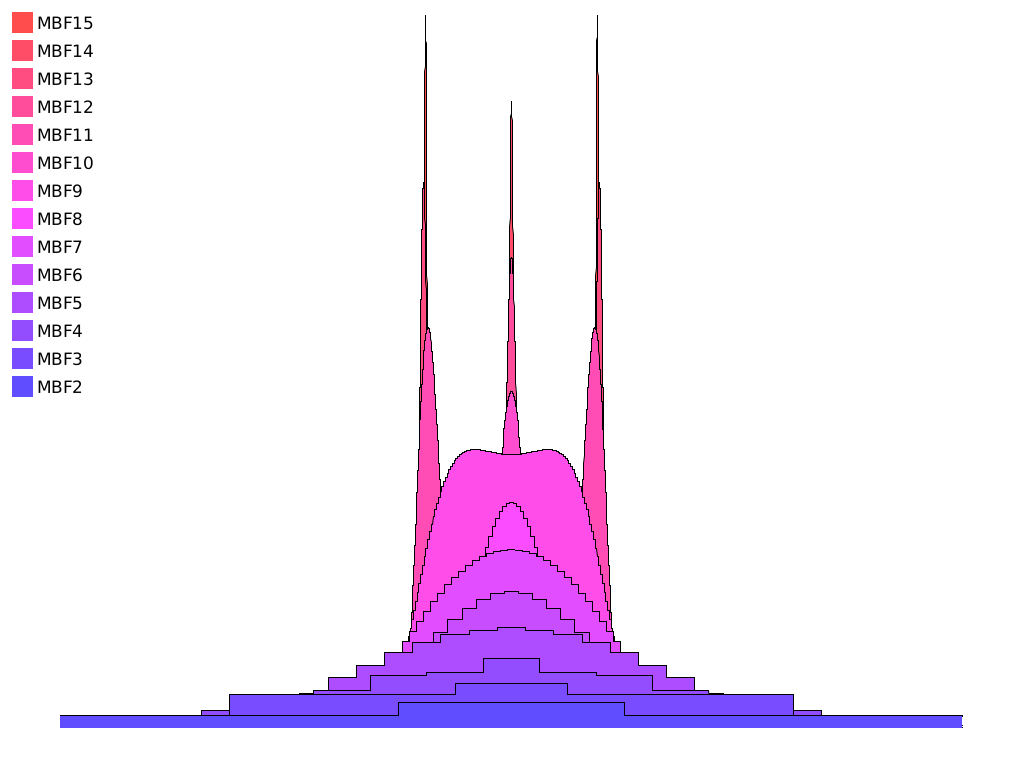}
        \caption{Real-valued y axis. Observe that at $M(9)$ the distribution begins to split into multiple lobes}
        \label{fig:layerSizeReal}
    \end{subfigure}
    \hfill
    \begin{subfigure}{0.48\linewidth}
        \centering
        \includegraphics[width=\linewidth]{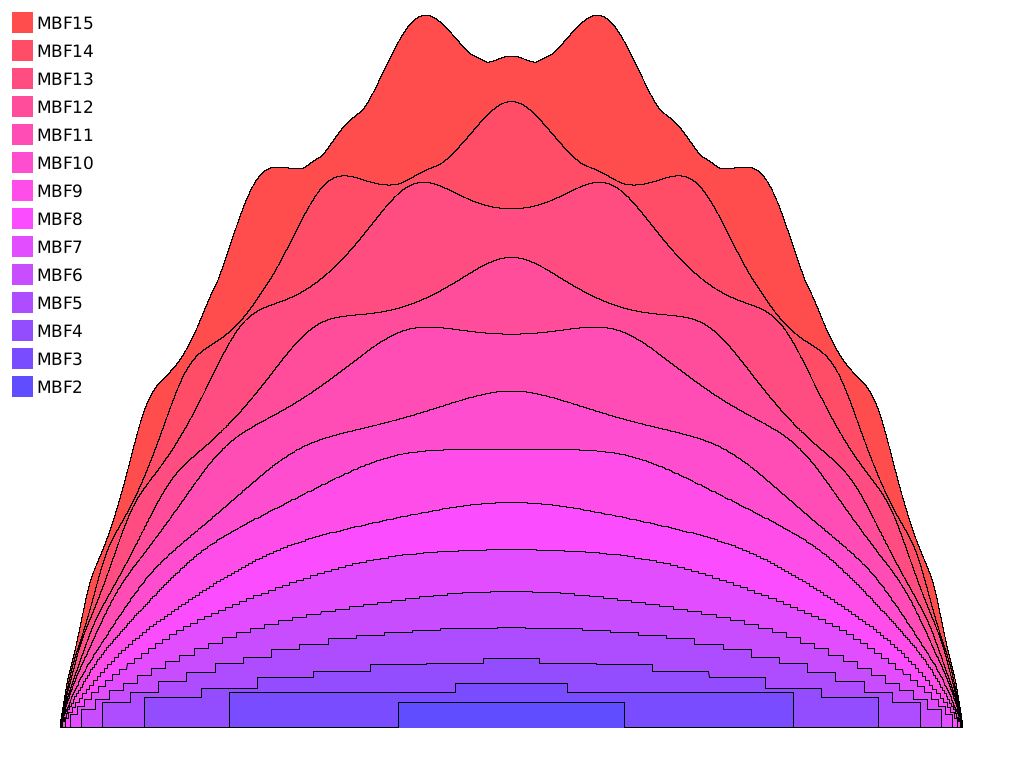}
        \caption{Logarithmic y axis. This depicts the harmonics within the layer counts of higher order $M(n)$ more clearly.}
        \label{fig:layerSizeLog}
    \end{subfigure}

    \caption{Both figures depict the relative number of MBFs in each layer. Horizontal axis is “progress from layer 0 to layer $2^n$”. The heights are scaled to the maximum value times an arbitrary constant per series chosen for visual clarity. The relative height between series has no meaning. Note that the curve for MBF9 is the uniform reference shown in Figure \ref{fig:mbf9DistributionComparison}.}
    \label{fig:layerSizeCombined}
\end{figure}

Starting from layer~$0$, we propagate the estimated branching ratios upward; starting from layer~$2^n$, we propagate downward. This continues until the two walks meet at the middle layer $2^{n-1}$. The relative discrepancy between the two independent estimates at this middle layer is given by
\[
\frac{S_{2^{n-1}}^{(0)}}{S_{2^{n-1}}^{(2^n)}} - 1,
\]
where $S_{2^{n-1}}^{(0)}$ denotes the estimate of cardinality of layer $2^{n-1}$ obtained from the walk beginning at layer~$0$, and $S_{2^{n-1}}^{(2^n)}$ denotes same cardinality estimate obtained from the walk beginning at layer~$2^n$. The relative sizes of these layers are illustrated by Figure \ref{fig:layerSizeCombined}.

The sum of the estimated layer cardinalities gives an estimate of $D(n)$.
By performing the experiment several times, we obtain an estimate of the standard error for each ratio,
and from there, a standard error for each layer cardinality estimate. We take the relative error of the most populated layer as an estimate of a standard error of the result. The results of these experiments and computational resources used are given in Table \ref{tbl:bestOtherDededkindEstimations}. 

\section{Results}
Through the use of three different methods, we were able to estimate Dedekind Numbers $D(10)$--$D(15)$ to two to four digit accuracy. Our most accurate results are listed in the introduction in Table~\ref{tab:dedekindNumbers}. Below we provide a breakdown for individual results by method. 
For $D(10)$ all methods we devised were applicable, and their results are compared in Table~\ref{tab:comparisonD10}. For $D(11)$ and $D(12)$, Pair matching was infeasible, so the remaining methods are compared in Table~\ref{tab:comparisonD11} and Table~\ref{tab:comparisonD12}. For $D(13)$, $D(14)$ and $D(15)$, only the weight layer branching method is feasible (the reference subset sampling method becomes impractical for $n \geq 13$; see Section~\ref{sect:wellKnownMBF}), and these results are given in Table~\ref{tbl:bestOtherDededkindEstimations}.

\begin{table}[H]
    \centering
    \begin{tabular}{l | c | c}
        \toprule
        Source & Value & S.E. \\
        \midrule
        Korshunov & $8.56 \times 10^{78}$ & n/a  \\
        Pair matching & $8.93782 \times 10^{78}$ & $2.234 \times 10^{75}$ \\
        Weight Layer Branching & $8.93708404 \times 10^{78}$ & $3.08476142 \times 10^{75}$ \\
        Reference Subsets & $8.93345 \times 10^{78}$ & $2.44 \times 10^{74} $ \\
        
        \hline
    \end{tabular}
    \caption{Comparison of D(10) values and standard errors from different experiments}
    \label{tab:comparisonD10}
\end{table}

\begin{table}[H]
    \centering
    \begin{tabular}{l | c | c}
        \toprule
        Source & Value & S.E. \\
        \midrule
        Korshunov & $2.56 \times 10^{144}$ & n/a  \\
        Weight Layer Branching & $3.6344\times 10^{144}$   & $1.68\times 10^{141}$ \\
        Reference Subsets & $ 3.6325 \times 10^{144} $ & $ 7.16 \times 10^{141} $ \\        
        \hline
    \end{tabular}
    \caption{Comparison of D(11) values and standard errors from different experiments}
    \label{tab:comparisonD11}
\end{table}

\begin{table}[H]
    \centering
    \begin{tabular}{l | c | c}
        \toprule
        Source & Value & S.E. \\
        \midrule
        Korshunov & $6.93 \times 10^{283}$ & n/a  \\
        Weight Layer Branching & $7.1491\times 10^{283}$   & $5.3\times 10^{280}$ \\
        Reference Subsets & $7.1911 \times 10^{283}$ & $3.25 \times 10^{281}$ \\        
        \hline
    \end{tabular}
    \caption{Comparison of D(12) values and standard errors from different experiments}
    \label{tab:comparisonD12}
\end{table}

\begin{table}[H]
\centering
\begin{tabular}{c|c|c|c|c}
\toprule
$n$ & Computed Value & S.E. & Runtime & Discrepancy \\
\midrule
10 & $8.93708404\times 10^{78}$    & $3.08476142\times 10^{75}$   & $6.1\,\mathrm{d}$         & $8.99\times 10^{-5}$ \\
11 & $3.63437475\times 10^{144}$   & $1.67617548\times 10^{141}$  & $13.3\,\mathrm{d}$        & $2.30\times 10^{-4}$ \\
12 & $7.14919263\times 10^{283}$   & $5.35341147\times 10^{280}$  & $28.7\,\mathrm{d}$        & $8.01\times 10^{-4}$ \\
13 & $6.03589738\times 10^{525}$   & $1.33138745\times 10^{523}$  & $36.6\,\mathrm{d}$        & $2.80\times 10^{-4}$ \\
14 & $5.58483742\times 10^{1043}$  & $4.06197320\times 10^{1041}$ & $110.7\,\mathrm{d}$ & $2.32\times 10^{-3}$ \\
15 & $3.80603932\times 10^{1953}$  & $5.29965888\times 10^{1951}$ & $260.21\,\mathrm{d}$        & $2.90\times 10^{-2}$ \\
\hline
\end{tabular}
\caption{Computational results of the Weight Layer Branching method for Dedekind numbers. S.E.\ denotes the standard error of the estimate. Runtime measured in core-days (\(\mathrm{d}\)) on an AMD Ryzen 7700x CPU.}
\label{tbl:bestOtherDededkindEstimations}
\end{table}

For $D(10)$, the pair-matching, reference-subset, and weight-layer methods all produce consistent estimates despite relying on very different ideas. Similar agreement is observed for $D(11)$ and $D(12)$, providing confidence that the reported values are close to the true Dedekind numbers. For $D(13)$--$D(15)$, only the weight-layer method is currently practical, yielding the first estimates for these dimensions together with explicit standard errors.

\section{Future Work}

\textbf{Metropolis-Hastings sampling.}
A Metropolis-Hastings variant of the current MCMC algorithm could be used to obtain a uniform stationary distribution directly, eliminating the need for degree correction. Since node degrees are already maintained during the walk, the additional cost is minimal.

\textbf{Analytical formula for expected Hamming distance.}
The expected Hamming distance between two independently and uniformly sampled MBFs is
\[
\mathbb{E}[d_H(f, g)] = 2 \sum_{w=0}^{n} \binom{n}{w} p_w (1 - p_w),
\]
where $p_w = P(f(x) = 1 \mid |x| = w)$ is the probability that a uniformly random MBF outputs 1 on an input of weight $w$. Since the values $p_w$ can be obtained directly from MCMC samples, this provides an alternative to the pairwise-distance experiments used to assess convergence.

\textbf{Larger reference subsets for $D(10)$.}
The current $D(10)$ estimate is limited by the low hit rate of the reference subset ($p_S \approx 1/1234$). Enlarging the subset while retaining an exact count would allow a larger fraction of MCMC samples to contribute to the estimate. Allowing limited variation in the layers adjacent to the middle layer appears to be a promising direction and could improve the accuracy of the estimate without increasing the amount of MCMC sampling.

\section*{Acknowledgments}
The authors thank the Paderborn Center for Parallel Computing (PC$^2$) for providing computational resources used in this work.

\section*{Declaration of generative AI and AI-assisted technologies in the manuscript preparation process.}
During the preparation of this work the authors used ChatGPT v.5.5 for final proofreading. After using this tool/service, the author(s) reviewed and edited the content as needed and are taking full responsibility for the content of the published article.

\appendix

\section{Buffer Matching Implementation Optimizations}
\label{app:bufferImpl}

\subsection{Bitset Optimization}
A key computational challenge with a naive two-loop iteration over all \texttt{left} and \texttt{right} MBFs is that it creates a very hot inner loop with an unpredictable branch. Especially since we cannot perform a 512-bit comparison in one instruction.

To make the comparison counting code branchless and more parallel, we perform a bitwise transpose of the MBFs in the \texttt{right} MBFs list, such that instead of having $|right|$ MBFs of 512 bits each, we instead have 512 bitsets of $|right|$ bits each.

Then, instead of comparing each element individually, we can do large Bitset ANDs between the bitsets corresponding to the \texttt{1} bits of the given \texttt{left} MBF. Then each \texttt{1} bit that remains in the resulting bitset corresponds to that \texttt{right} MBF being a valid match for this \texttt{left} MBF. To count, we simply popcount the resulting bitset, which can be implemented efficiently using CPU intrinsics.

We can improve this further by exploiting the fact that we're working with MBFs, and thus there are implication relations between the bits. We only need to check the bits corresponding to the \texttt{1} bits in the antichain representation of the chosen \texttt{left} MBF. In our measurements, we found the average number of \texttt{1} bits in the antichain representation for 9 variables is around 51.857 (see Table \ref{tab:averageMinCuts}).

\begin{table}[H]
    \centering
    \begin{tabular}{c|c|c|c}
        \textbf{Vars} & \textbf{Avg size} & \textbf{S.E.} & \textbf{Exact total} \\
        1 & 0.666666 & 0 & 2/3 \\
        2 & 1.0 & 0 & 6/6 \\
        3 & 1.6 & 0 & 32/20 \\
        4 & 2.70238 & 0 & 454/168 \\
        5 & 4.68434 & 0 & 35512/7581 \\
        6 & 8.50554 & 0 & 66584412/7828354 \\
        7 & 15.1433 & 0 & 36566354210304/2414682040998 \\
        8 & 29.7951717 & 0.0000222 &  \\
        9 & 51.8571551 & 0.0000503 &  \\
        10 & 112.100685 & 0.0000801 &  \\
        11 & 202.605450 & 0.000345 &  \\
        12 & 429.118260 & 0.000255 &  \\
        13 & 801.220772 & 0.000558 &  \\
        14 & 1387.63473 & 0.00141 &  \\
        15 & 2294.66560 & 0.767 &  \\
    \end{tabular}
    \caption{Average antichain sizes by number of variables. The first 7 values were computed exactly from the dataset of all MBFs\cite{vanhirtum2023computation}. Note that the exact totals correspond to \cite{oeisA118077}. Terms 8 and higher were computed using the MCMC method described in Section \ref{sect:walks}.}
    \label{tab:averageMinCuts}
\end{table}

This reduces the number of AND operations from 256 to 51.86 on average.

\subsection{AVX2 SIMD Optimizations}
\label{sect:simdBufferComparison}

We performed a parameter study to determine the optimal stride of AND operations. Table \ref{tab:simdStrideBenchmarks} shows results for various SIMD block sizes in bits. The optimal size is 2048 bits (8 AVX2 registers), corresponding to 4 cache lines. Larger register counts incur stack spills that negate the benefit.

\begin{table}[H]
    \centering
    \begin{tabular}{c|c}
        \textbf{SIMD Block Size (bits)} & \textbf{Time Taken (s)} \\
        256  & 35.5012 \\
        512  & 19.0837 \\
        1024 & 18.0097 \\
        1536 & 18.5889 \\
        2048 & 16.7531 \\
        3072 & 21.1548 \\
        4096 & 21.3835 \\
        8192 & 29.3380 \\
    \end{tabular}
    \caption{Performance of SIMD-based buffer matching for different SIMD block sizes. All variants found the exact same number of valid combinations: 29938673/274877906944=0.0108916\%.}
    \label{tab:simdStrideBenchmarks}
\end{table}

\subsection{Speedup through Comparison Pre-filtering}
Since $p_{match,9} \approx 10^{-4}$, the vast majority of pairs will not form a valid MBF10. Pre-filtering partitions the MBFs into categories based on a subset of sampled bits; any category pair that is already incompatible can be dismissed without examining the individual MBFs.

\subsection{Initial Table-based Approach}
The initial approach selected 16 sample bits from the two middle layers and assigned each \texttt{right} MBF to one of $2^{16}=65536$ categories. For a given \texttt{left} MBF, any category whose sample bits are not a subset of $sample_{left}$ is skipped. Ideally, the fraction of comparisons that must be examined is approximated by

\begin{equation}
    \label{eq:fractionCheck}
    fraction\_check = \frac{\bigsum{n}{i = 0}{i! 2^i}}{2^{2n}}
\end{equation}

Because MBF9 bits in the two middle layers (L4 and L5) are strongly biased (L4: $\approx$21.5\% ones; L5: $\approx$78.5\% ones), we AND together two random L4 bits and OR together two random L5 bits to produce 8 sample bits from each layer. The resulting filtering fractions are shown in Table~\ref{tab:filterFractions}.

\begin{table}[H]
    \centering
    \begin{tabular}{c|l}
        1 bit & 12.4666\% \\
        2 bits & 9.19239\% \\
        3 bits & 10.2990\% \\
        4 bits & 12.4734\% \\
        5 bits & 16.0678\% \\
    \end{tabular}
    \caption{Filtering percentages for various numbers of bits ANDed (for L4) or ORed (for L5) on a dataset of $(2^{28})^2$ \texttt{left} and \texttt{right} MBFs.}
    \label{tab:filterFractions}
\end{table}

Filtering below 9.2\% could not be achieved with the table method, because the extreme (000\ldots0 and 111\ldots1) categories are overrepresented due to the bimodal weight distribution visible in Figure~\ref{fig:mbf9DistributionComparison}.

\subsection{Filtering Binary Tree Structure}
\label{sect:treeImpl}

\begin{figure}[H]
    \centering
    \includegraphics[width=0.6\linewidth]{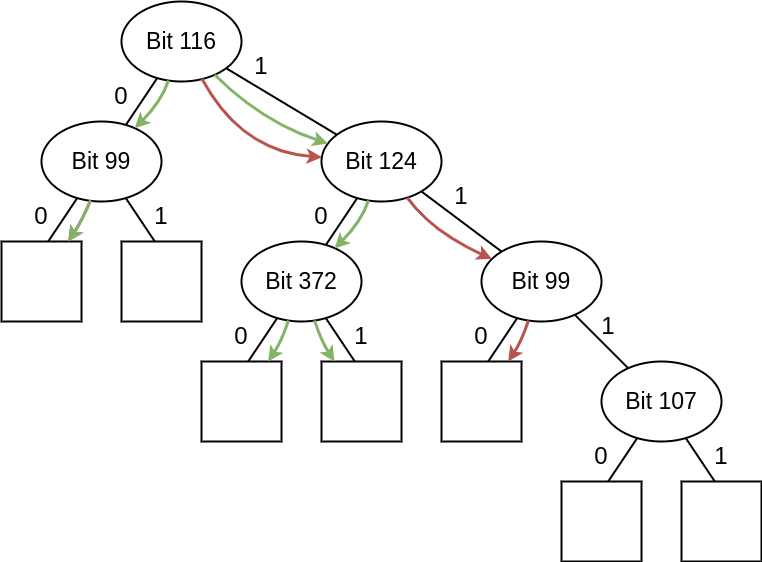}
    \caption{Example of the tree filter. Round nodes represent a sample-bit query; square leaf nodes are the buckets MBFs are placed in. \texttt{left} MBF traversal shown in green, \texttt{right} MBF traversal in red.}
    \label{fig:treeFilterExample}
\end{figure}

Instead of a flat table, we structure the MBF categories as a binary tree whose interior nodes each split on a chosen MBF bit. \texttt{right} MBFs follow a single branch based on the bit value. \texttt{left} MBFs follow all branches whose bit value is $\leq$ the left MBF's bit, splitting at every node where the left bit is 1. A right MBF whose leaf cannot be reached by a given left MBF is filtered out without a full comparison.

\subsection{Constructing the Tree}

\subsubsection{Splitting Criterion}
At each node, we choose the bit that maximises the number of eliminated comparisons over a sample of 255 left and right MBFs:

\begin{equation}
    \label{eq:splitCriterion}
    \underset{bit}{\operatorname{argmax}}\ \#\text{RMBF}[bit=1]\cdot\#\text{LMBF}[bit=0]
\end{equation}

Splitting stops when the expected savings fall below a threshold of $131{,}072{,}000$ (\texttt{MIN\_SPLIT\_COUNT}), chosen empirically to match the SIMD block size at which the bitset comparison (Appendix~\ref{sect:simdBufferComparison}) becomes efficient. At the 10M block size, the tree eliminates 92.46\% of comparisons (filtering percentage 7.54\%); at the 10k block size, only 44.96\% of comparisons are eliminated (filtering percentage 55.04\%).


\bibliographystyle{elsarticle-num} 
\bibliography{refs.bib}

\end{document}